\documentclass[11pt,twoside]{amsart}

\usepackage[T1]{fontenc} \usepackage[latin1]{inputenc}

\usepackage{amsmath}
\usepackage{amsfonts}
\usepackage{amssymb}
\usepackage{amsthm}
\usepackage{tikz}
\usepackage[all,cmtip,color,arc,curve]{xy}

\theoremstyle{definition}
\newtheorem{defi}{Definition}[section]

\theoremstyle{plain}
\newtheorem{thm}[defi]{Theorem}
\newtheorem{lem}[defi]{Lemma}
\newtheorem{prop}[defi]{Proposition}

\newtheorem{conj}[defi]{Conjecture}

\newcommand*{\ch}{\ensuremath{\operatorname{ch}}}
\newcommand*{\Coh}{\ensuremath{\operatorname{Coh}}}
\newcommand*{\into}{\ensuremath{\hookrightarrow}}
\newcommand*{\onto}{\ensuremath{\twoheadrightarrow}}

\newcommand*{\Hom}{\ensuremath{\operatorname{Hom}}}
\newcommand*{\Ext}{\ensuremath{\operatorname{Ext}}}

\newcommand*{\calA}{\ensuremath{\mathcal{A}}}
\newcommand*{\calB}{\ensuremath{\mathcal{B}}}
\newcommand*{\calC}{\ensuremath{\mathcal{C}}}

\newcommand*{\calF}{\ensuremath{\mathcal{F}}}
\newcommand*{\calO}{\ensuremath{\mathcal{O}}}

\newcommand*{\calT}{\ensuremath{\mathcal{T}}}
\newcommand*{\bbA}{\ensuremath{\mathbb{A}}}
\newcommand*{\bbC}{\ensuremath{\mathbb{C}}}
\newcommand*{\bbP}{\ensuremath{\mathbb{P}}}
\newcommand*{\bbQ}{\ensuremath{\mathbb{Q}}}
\newcommand*{\bbR}{\ensuremath{\mathbb{R}}}
\newcommand*{\bbZ}{\ensuremath{\mathbb{Z}}}

\begin{document}

\title[Bogomolov-Gieseker inequality for the quadric threefold]{A generalized
Bogomolov-Gieseker inequality for the smooth quadric threefold}

\author{Benjamin Schmidt}
\address{Department of Mathematics, The Ohio State University, 231 W 18th
Avenue, Columbus, OH 43210-1174, USA}
\email{schmidt.707@osu.edu}
\urladdr{https://people.math.osu.edu/schmidt.707/}



\begin{abstract}
We prove a generalized Bogomolov-Gieseker inequality as conjectured by Bayer,
Macr\`i and Toda for the smooth quadric threefold. This implies the existence of
a family of Bridgeland stability conditions.
\end{abstract}

\maketitle

\section{Introduction}

The classical notion of $\mu$-stability has been explored for a long time to
study vector bundles and their moduli spaces. One important direction of study
is the birational geometry of a given moduli space. Historically, an approach
for obtaining divisorial contractions or flips was varying the
polarization of the variety and therefore varying the GIT problem. However, this
does not provide enough flexibility. For example, if the Picard group is $\bbZ$,
there is no possible variation.

Inspired by the study of Dirichlet branes in string theory by Douglas (see
\cite{Dou00, Dou01, Dou02}), the notion of Bridgeland stability was introduced
in \cite{Bri07}. Instead of defining stability in the category of coherent
sheaves, one uses other abelian categories inside the bounded derived category
of coherent sheaves. Bridgeland shows that the set of all these stability
conditions forms a complex manifold. This leads to plenty of room to vary a
given stability condition even if the Picard rank is $1$.

While this notion provides many of the desired properties, constructing such
Bridgeland stability conditions has turned out to be a serious issue. A large
family was constructed in the case of K3 surfaces in \cite{Bri08}. Arcara
and Bertram generalized this construction to any smooth complex
projective surface in \cite{AB13}. Examples of successful applications are found
in the birational geometry of Hilbert schemes of points on smooth projective
surfaces (see for example \cite{ABCH13, BM13, MM13, YY14}). Toda shows that the
minimal model program on any smooth projective surface is realized as a
variation of moduli spaces of Bridgeland stable objects in \cite{Tod12}.

The case of threefolds seems to be more complicated. The work of Bridgeland was
motivated by the case of Calabi-Yau threefolds occurring in string theory. So
far no Bridgeland stability condition has been constructed on a single
Calabi-Yau threefold. A promising approach for all smooth projective threefolds
is due to Bayer, Macr\`i and Toda in \cite{BMT14}. It was confirmed to work for $\bbP^3$
in \cite{Mac12} and for principally polarized abelian threefolds of Picard rank
one in \cite{MP13a, MP13b}. By mimicking the construction for surfaces, Bayer,
Macr\`i and Toda obtain the notion of tilt-stability on an abelian category
$\calB^{\omega, B}$ in the bounded derived category of coherent sheaves for any
$\bbR$-divisor $B$ and any ample $\bbR$-divisor $\omega$. The slope function is
given by
$$\nu_{\omega, B} := \frac{\omega \ch^B_2 - \frac{\omega^3}{2} \ch^B_0}{\omega^2
\ch^B_1},$$
where $\ch^B = e^{-B} \ch$. Unlike in the case of surfaces this provides no
Bridgeland stability condition. They conjecture a generalized Bogomolov-Gieseker
inequality on third Chern classes for tilt-stable objects $E \in \calB^{\omega,
B}$ which satisfy $\nu_{\omega, B}(E) = 0$ given by $$\ch^B_3(E) \leq
\frac{\omega^2}{6} \ch^B_1(E).$$
This inequality turns out to be the missing ingredient for the
construction of Bridgeland stability conditions. Interestingly,
there are other applications of this inequality besides the construction of
Bridgeland stability conditions. One of the most interesting consequences is
Fujita's conjecture (see \cite{BBMT11}). Macr\`i was able to prove the
inequality in the case of $\bbP^3$ in \cite{Mac12}, while Maciocia and Piyaratne
managed to show it for principally polarized abelian threefolds of Picard rank
one in \cite{MP13a, MP13b}. The main result of this article is the following.

\begin{thm}(See Theorem \ref{thm-main})
The generalized Bogomolov-Gieseker inequality is true for the smooth quadric
threefold $Q$. In particular, there is a large family of Bridgeland stability
conditions on $Q$.
\end{thm}

The proof is based on calculations with a strong full exceptional collection in
$D^b(Q)$ that exists due to \cite{Kap88}. We break it down to a technical
lemma from \cite{BMT14} (see Proposition \ref{basic-lemma}).

The paper is organized as follows. Basics on stability and the
construction of \cite{BMT14} are explained in Section \ref{sec:construction}.
In Section \ref{sec:quadric} some facts about the smooth quadric threefold are
being recalled. Finally, Section \ref{sec:main} deals with the proof of the main
theorem.

\subsection*{Notation}

By $X$ we denote a smooth projective threefold over the complex numbers. Its
bounded derived category of coherent sheaves is called $D^b(X)$. Let $Q$ be
the smooth quadric threefold in $\bbP^4$ over the complex numbers defined
by the equation $x_0^2 + x_1 x_2 + x_3 x_4 = 0$. 

\subsection*{Acknowledgements}
I would like to thank Emanuele Macr\`i for reading preliminary versions of this
article and for many useful discussions. I also appreciate useful advice from
Arend Bayer on a previous version of this manuscript. I thank the referee for
the detailed reading of this article. The research was partially supported by
NSF grants DMS-1160466 and DMS-1302730 (PI Emanuele Macr\`i).

%The research was supported by NSF grant
%DMS-1160466 (PI Emanuele Macr\`i).

\section{Construction of Stability Conditions}
\label{sec:construction}

Let us recall some definitions concerning stability. The central part of the
theory is the notion of Bridgeland stability conditions that was introduced in
\cite{Bri07}. Let $H:=\{re^{i\pi\varphi} : r>0, \varphi \in (0,1]\}$ be the
upper half plane plus the negative real line. A Bridgeland stability condition
on $D^b(X)$ is a pair $(Z,\calA)$, where $\calA$ is the heart of a bounded
t-structure and $Z: K_0(X) = K_0(\calA) \to \bbC$ is a homomorphism such that
$Z(\calA \backslash \{ 0 \}) \subset H$ holds plus a technical property.

The inclusion $Z( \calA \backslash \{ 0 \}) \subset H$ turns out to be the
crucial point for threefolds. Note that for any smooth projective variety of
dimension bigger than or equal to two, there is no Bridgeland stability
condition factoring through the Chern character for $\calA = \Coh(X)$ due to \cite[Lemma 2.7]{Tod09}.

In order to construct such stability conditions on a smooth projective
threefold $X$, Bayer, Macr\`i and Toda proposed a construction in \cite{BMT14}.
We will review it. Let $B$ be any $\bbR$-divisor. Then the twisted Chern
character $\ch^B$ is defined to be $e^{-B} \ch$. In more detail, we have
\begin{align*}
\ch^B_0 &= \ch_0, \\
\ch^B_1 &= \ch_1 - B \ch_0 ,\\
\ch^B_2 &= \ch_2 - B \ch_1 + \frac{B^2}{2} \ch_0, \\
\ch^B_3 &= \ch_3 - B \ch_2 + \frac{B^2}{2} \ch_1 - \frac{B^3}{6} \ch_0.
\end{align*}

The category of coherent sheaves $\Coh(X)$ is the heart of a
bounded t-structure on $D^b(X)$. Let $\omega$ be an ample $\bbR$-divisor. Then
we can define a twisted version of the standard slope stability function on
$\Coh(X)$ by
$$\mu_{\omega,B} := \frac{\omega^2 \ch^B_1}{\omega^3 \ch^B_0},$$
where dividing by $0$ is interpreted as $+\infty$.
%A sheaf $E \in \coh(X)$ is called $\mu_{\omega,B}-$ or slope-(semi)stable if
% for all subsheaves $F \subset E$ the inequality $\mu_{\omega,B}(F) < (\leq)
% \mu_{\omega,B}(E/F)$ holds.
The process of tilting is used to obtain a new heart of a bounded t-structure.
For more information on this general theory we refer to \cite{HRS96}. A torsion
pair is defined by
\begin{align*}
\calT_{\omega, B} &= \{E \in \Coh(X) : \text{any quotient $E \onto G$
satisfies $\mu_{\omega,B}(G) > 0$} \}, \\
\calF_{\omega, B} &=  \{E \in \Coh(X) : \text{any subsheaf $F \subset E$
satisfies $\mu_{\omega,B}(F) \leq 0$} \}.
\end{align*}
A new heart of a bounded t-structure is defined by the extension
closure $\calB^{\omega, B}(X) := \langle \calF_{\omega, B}[1],
\calT_{\omega, B} \rangle$. A new slope function is defined by
$$\nu_{\omega, B} := \frac{\omega \ch^B_2 - \frac{\omega^3}{2}
\ch^B_0}{\omega^2 \ch^B_1},$$
where dividing by $0$ is again interpreted as $+\infty$. Note that in
regard to \cite{BMT14} this slope has been modified by switching $\omega$ with
$\sqrt{3} \omega$. We prefer this point of view because it will slightly
simplify a few computations. On smooth projective surfaces the map $-\ch^B_2 +
\frac{\omega^2}{2} \ch^B_0 + i \omega \ch^B_1$ from $\calB^{\omega, B}(X)$ to
$\bbC$ is already a Bridgeland stability function (see \cite{Bri08, AB13}).
However, on threefolds this is not enough. For example skyscraper sheaves are
still mapped to the origin. Therefore, Bayer, Macr\`i and Toda propose another
analogous tilt via
\begin{align*}
\calT'_{\omega, B} &= \{E \in \calB^{\omega, B}(X) : \text{any quotient $E \onto
G$ satisfies $\nu_{\omega,B}(G) > 0$} \}, \\
\calF'_{\omega, B} &=  \{E \in \calB^{\omega, B}(X) : \text{any subobject $F
\into E$ satisfies $\nu_{\omega,B}(F) \leq 0$} \}
\end{align*}
and setting $\calA^{\omega,B} := \langle \calF'_{\omega, B}[1],
\calT'_{\omega, B} \rangle $. Finally, they define for any
$s>0$ functions by
\begin{align*}
Z_{\omega,B,s} &:= (-\ch^B_3 + s\omega^2\ch^B_1) + i (\omega \ch^B_2 -
\frac{\omega^3}{2}\ch^B_0), \\
\lambda_{\omega,B,s} &:= -\frac{\Re(Z_{\omega,B,s})}{\Im(Z_{\omega,B,s})}.
\end{align*}
The function $\lambda_{\omega,B,s}$ is called the slope of $Z_{\omega,B,s}$.

\begin{defi}
An object $E \in \calB^{\omega, B}$ is called \textit{$\nu_{\omega,
B}$-(semi)stable} (or \textit{tilt-(semi)stable}) if for any exact sequence $0
\to F \to E \to G \to 0$ the inequality $\nu_{\omega,B}(F) < (\leq) \nu_{\omega,B}(G)$ holds.
\end{defi}

The following theorem motivates the whole construction. 

\begin{thm}[{\cite[Corollary 5.2.4]{BMT14}}]
\label{BMT-main}
Let $X$ be a smooth projective threefold over the complex numbers, $\omega$ an
ample divisor, $B$ any divisor and $s>0$. Then $(Z_{\omega, B, s}
,\calA^{\omega, B})$ is a Bridgeland stability condition if and only if for any
$\nu_{\omega,B}$-stable object $E \in \calB^{\omega,B}$ with
$\nu_{\omega,B}(E)=0$ the inequality
\begin{align}
\label{main-ineq}
\ch^B_3(E) < s \omega^2 \ch^B_1 (E)
\end{align}
holds.
\end{thm}

The inequality \eqref{main-ineq} in the theorem is exactly expressing the fact
that $Z_{\omega, B, s}$ is not mapping on the non-negative real line $\bbR_{\geq
0}$. Bayer, Macr\`i and Toda hope that \eqref{main-ineq} holds for $s=3/2$. They
even conjecture a stronger inequality.

\begin{conj}[{\cite[Conjecture 1.3.1]{BMT14}}]\label{conj-main}
Inequality \eqref{main-ineq} holds for all $s>\frac{1}{6}.$
\end{conj}
\section{Quadric Threefold}
\label{sec:quadric}

In order to prove Conjecture \ref{conj-main} for the smooth quadric threefold
$Q$, we need to recall some facts about its bounded derived category of coherent
sheaves $D^b(Q)$. In the following, we view $Q$ as being cut out by the equation
$x_0^2 + x_1 x_2 + x_3 x_4 = 0$ in $\bbP^4$.

Since the open subvariety of $Q$ defined by $x_1 \neq 0$ is isomorphic to
$\bbA^3$, the Picard group of $Q$ is isomorphic to $\bbZ$ and is
generated by a very ample line bundle $\calO(H)$. Moreover, the equality
$H^3 = 2$ holds because a general line in $\bbP^4$ intersects $Q$ in two points.

Let us recall exceptional collections.

\begin{defi}
  A \textit{strong exceptional collection} is a sequence $E_1, \ldots, E_r$ of
  objects in $D^b(X)$ such that $\Ext^i(E_l, E_j)
  = 0$ for all $l,j$ and $i\neq 0$, $\Hom(E_j, E_j)=\bbC$ and $\Hom(E_l, E_j) =
  0$ for all $l>j$. Moreover, it is called \textit{full} if $E_1, \ldots, E_r$
  generates $D^b(X)$ via shifts and extensions.
\end{defi}

On $Q$ line bundles are not enough to obtain a full strong exceptional
collection. Therefore, we need to introduce the spinor bundle $S$. We refer to
\cite{Ott88} for a more detailed treatment. The spinor bundle is defined via an
exact sequence
$$0 \to \calO_{\bbP^4}(-1)^{\oplus 4} \to \calO_{\bbP^4}^{\oplus 4} \to \iota_*
S \to 0$$
where $\iota: Q \into \bbP^4$ is the inclusion and the first map is given by a
matrix $M$ such that $M^2 = (x_0^2 + x_1 x_2 + x_3 x_4) I_4$ for the identity
$4\times 4$ matrix $I_4$. Restricting the second morphism to $Q$ leads to
\begin{align}
\label{spinor-seq}
0 \to S(-1) \to \calO_Q^{\oplus 4} \to S \to 0.
\end{align}
Due to Kapranov (see \cite{Kap88})
$$\calO(-1), S(-1), \calO, \calO(1)$$
is a strong full exceptional collection on $D^b(Q)$.

Explicit computations lead to a resolution of the skyscraper sheaf $k(x)$ given
by
\begin{align}
\label{eq:resolution}
0 \to \calO(-1) \to S(-1)^{\oplus 2} \to \calO^{\oplus 4} \to
\calO(1) \to k(x) \to 0
\end{align}
for any $x \in Q$.

\section{Main Result}
\label{sec:main}

The main result of this article is the following.

\begin{thm}\label{thm-main}
Conjecture \ref{conj-main} holds for the smooth projective threefold $Q$,
i.e., for any $\nu_{\omega,B}$-stable object $E \in \calB^{\omega,B}$ with
$\nu_{\omega,B}(E)=0$ the inequality
\begin{align*}
\ch^B_3(E) \leq \frac{\omega^2}{6} \ch^B_1 (E)
\end{align*}
holds.
\end{thm}

There are $\alpha \in \bbR_{> 0}$ and $\beta \in \bbR$ such that $\omega =
\alpha H$ and $B = \beta H$. Therefore, we will replace $B$ by $\beta$ and $\omega$ by
$\alpha$ in the notation of slope functions and categories. Due to Proposition
2.7 and Lemma 3.2 in \cite{Mac12} it suffices to prove the statement for $\alpha
< \frac{1}{3}$ and $\beta \in [-\frac{1}{2}, 0]$.

%Moreover, we will assume that
%$\beta \neq -\alpha$ and deal with this special case later.

The following technical proposition provides the basis of the proof.

\begin{prop}[{\cite[Lemma 8.1.1]{BMT14}}]
\label{basic-lemma}
Let $\calC \subset D^b(X)$ be the heart of a bounded t-structure with the
following properties.
\begin{enumerate}
  \item There exists $\phi_0 \in (0,1)$ and $s_0 \in \bbQ$ such that
  $$ Z_{\alpha,\beta,s_0} (\calC) \subset \{r e^{\pi \phi i} : r \geq 0,
  \phi_0 \leq \phi \leq \phi_0 +1 \}.$$
  \item The inclusion $\calC \subset \langle \calA_{\alpha,\beta},
  \calA_{\alpha,\beta}[1] \rangle$ holds.
  \item For all points $x \in X$ we have $k(x) \in \calC$ and for all proper
  subobjects $C \into k(x)$ in $\calC$ the inequality $\Im
  Z_{\alpha,\beta,s_0}(C) > 0$ holds.
\end{enumerate}
Then the pair $(Z_{\alpha,\beta,s}, \calA_{\alpha, \beta})$ is a stability
condition on $D^b(X)$ for all $s>s_0$.
\end{prop}

Due to \cite{Bon90} a full strong exceptional collection induces an equivalence
between $D^b(X)$ and the bounded derived category of finitely generated
modules over some finite dimensional algebra $A$. In the special case of the
smooth quadric $Q$, we get the heart of a bounded t-structure by setting
$$\calC := \langle \calO(-1)[3], S(-1)[2], \calO [1], \calO(1) \rangle.$$
Moreover, $\calC$ is isomorphic to the category of finitely generated
modules over some finite dimensional algebra $A$ and $\calO(-1)[3], S(-1)[2],
\calO [1], \calO(1)$ are the simple objects.

We will show that the conditions of the lemma are fulfilled for this $\calC$ and
$s_0 = \frac{1}{6}$. In order to do that, a computation of the values for the
different slope-functions is necessary. By using \eqref{spinor-seq} we can
obtain the following lemma.

\begin{lem}
\label{computation}
For all $n \in \mathbb{N}$ we have
$$\ch^{\beta} (\calO(n)) = 1 + (n-\beta)H + (n-\beta)^2 \frac{H^2}{2} +
\frac{1}{3}(n-\beta)^3.$$
The chern character of $S(-1)$ is given by
$$\ch^{\beta} (S(-1)) = 2 - (2\beta + 1)H + \beta(\beta + 1) H^2 +
\frac{1}{6} - \beta^2 - \frac{2}{3} \beta^3.$$
We have the following $\mu$-slopes
\begin{alignat*}{2}
&\mu_{\alpha,\beta}(\calO(1)) = \frac{1-\beta}{\alpha}, \
&&\mu_{\alpha,\beta}(\calO) = -\frac{\beta}{\alpha}, \\
&\mu_{\alpha,\beta}(\calO(-1)) = -\frac{\beta+1}{\alpha}, \
&&\mu_{\alpha,\beta}(S(-1)) = -\frac{2\beta + 1}{2 \alpha}.
\end{alignat*}
The $\nu$-slopes for the same sheaves are given by
\begin{alignat*}{2}
&\nu_{\alpha,\beta}(\calO(1)) = \frac{(1-\beta)^2 - \alpha^2}{2\alpha(1-\beta)},
\ &&\nu_{\alpha,\beta}(\calO) = \frac{\alpha^2 - \beta^2}{2\alpha \beta}, \\
&\nu_{\alpha,\beta}(\calO(-1)) = \frac{\alpha^2 - (1+\beta)^2}{2\alpha
(1+\beta)}, \ &&\nu_{\alpha,\beta}(S(-1)) = \frac{\alpha^2 - \beta (\beta +
1)}{\alpha (2\beta + 1)}.
\end{alignat*}
Finally, the $Z$ values can be computed as
\begin{align*}
&Z_{\alpha,\beta,\frac{1}{6}}(\calO(1)) = \frac{1}{3}((1-\beta)^2 - \alpha^2)(\beta
-1 + 3 i \alpha), \\
&Z_{\alpha,\beta,\frac{1}{6}}(\calO) = \frac{1}{3}(\beta^2 - \alpha^2)(\beta + 3
i \alpha), \\
&Z_{\alpha,\beta,\frac{1}{6}}(\calO(-1)) = \frac{1}{3}((1+\beta)^2 -
\alpha^2)(\beta + 1 + 3i\alpha), \\
&Z_{\alpha,\beta,\frac{1}{6}}(S(-1)) = \frac{1}{6}(2\beta + 1)(2\beta^2 + 2\beta -
1 - 2 \alpha^2) + 2 i \alpha (\beta^2 + \beta - \alpha^2).
\end{align*}
\end{lem}

At this point we can prove the first assumption in Proposition
\ref{basic-lemma}.

\begin{lem}
There exists $\phi_0 \in (0,1)$ such that
$$Z_{\alpha,\beta,\frac{1}{6}} (\calC) \subset \{r e^{\pi \phi i} : r \geq 0,
\phi_0 \leq \phi \leq \phi_0 + 1 \}.$$
\begin{proof}
It suffices to show that the $4$ generators of $\calC$ are contained in some
half plane of $\bbC$. There are two different cases to deal with. Lemma
\ref{computation} shows that the half plane of points with negative real part
works if $|\beta| \leq |\alpha|$, while the half plane left of the line through
$0$ and $Z_{\alpha,\beta,\frac{1}{6}} (\calO[1])$ works in the case $|\beta| >
|\alpha|$. The following figure shows the $Z_{\alpha,\beta,\frac{1}{6}}$ values.

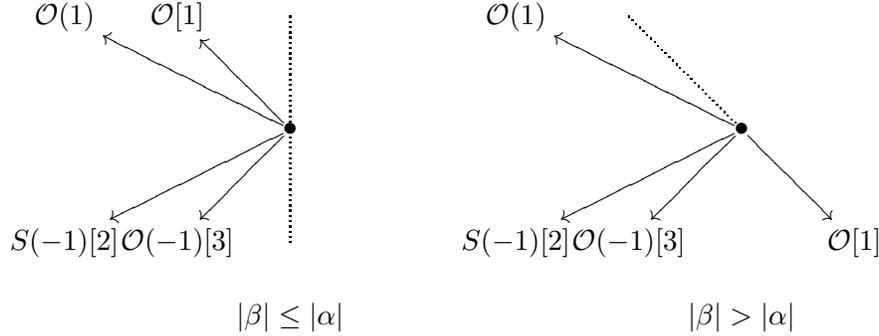
\centerline{
\xygraph{
!{<0cm,0cm>;<1cm,0cm>:<0cm,1cm>::}
%Vertices 1
!{(0,0) }*+{\calO(1)}="1"
!{(1.5,0) }*+{\calO[1]}="2"
!{(3,0) }*{}="3"
!{(3,-1.5) }*{\bullet}="4"
!{(0,-3) }*+{S(-1)[2]}="5"
!{(1.5,-3) }*+{\calO(-1)[3]}="6"
!{(3,-3) }*{}="7"
!{(3,-4) }*{|\beta| \leq |\alpha|}="text"
%Vertices 2
!{(6,0) }*+{\calO(1)}="1'"
!{(10.5,-3) }*+{\calO[1]}="2'"
!{(7.5,0) }*{}="3'"
!{(9,-1.5) }*{\bullet}="4'"
!{(6,-3) }*+{S(-1)[2]}="5'"
!{(7.5,-3) }*+{\calO(-1)[3]}="6'"
!{(9,-4) }*{|\beta| > |\alpha|}="text'"
%Edges 1
"3"-@{.}"7"
"4"-@{>}"1"
"4"-@{>}"2"
"4"-@{>}"5"
"4"-@{>}"6"
%Edges 2
"3'"-@{.}"4'"
"4'"-@{>}"1'"
"4'"-@{>}"2'"
"4'"-@{>}"5'"
"4'"-@{>}"6'"
}}
\end{proof}
\end{lem}

Before we can show assumption (ii) in Proposition \ref{basic-lemma}, we need to
deal with continuity issues for tilt-stability. For any $E \in
\calB^{\alpha,\beta}$ we denote the minimum of all $\nu_{\alpha, \beta}(G)$ for
quotients $E \onto G$ by $\nu^{\min}_{\alpha, \beta}(E)$.

\begin{lem}
\label{lem:open}
Let $E \onto N$ be an epimorphism in the category $\calB^{\alpha_0,\beta_0}$
where $N$ is the semistable quotient in the Harder-Narasimhan filtration. Assume
additionally that $E$ has no subobject with $\nu_{\alpha_0,\beta_0} = \infty$.
Then there is an open subset $U$ around the point $(\alpha_0,\beta_0)$ such that
the following holds.
\begin{enumerate}
  \item The inequality $\nu^{\min}_{\alpha, \beta}(E) \leq \nu_{\alpha,
\beta}(N)$ holds for all $(\alpha, \beta) \in U$.
  \item If $N_1, \ldots, N_l$ are the stable factors of $N$, then we obtain the
  inequality $\nu^{\min}_{\alpha, \beta}(E) \geq \min \nu_{\alpha,\beta}(N_i)$
  for all $(\alpha, \beta) \in U$.
\end{enumerate}
\begin{proof}
By definition we have $\nu^{\min}_{\alpha_0, \beta_0}(E) = \nu_{\alpha_0,
\beta_0}(N)$. Each semistable factor in the Harder-Narasimhan filtration of $E$
has a Jordan-H\"older filtration by stable factors. Since none of these stable
factors has $\nu_{\alpha_0,\beta_0} = \infty$, we can use openness of stability
(\cite[Corollary 3.3.3]{BMT14}) to show that all these stable factors are in the
category $\calB^{\alpha,\beta}$ in a small open neighborhood $U$ of
$(\alpha_0,\beta_0)$. But that means $E$, $N$ and the kernel of $E \onto N$
are in $\calB^{\alpha,\beta}$ for all ${\alpha, \beta} \in U$. Therefore, we
have $E \onto N$ in $\calB^{\alpha,\beta}$ for all ${\alpha, \beta} \in U$.
But that implies $\nu^{\min}_{\alpha, \beta}(E) \leq \nu_{\alpha, \beta}(N)$ for all
$(\alpha, \beta) \in U$.

We shrink $U$ such that the slopes of the $N_i$ are smaller than the
slopes of all the other stable factors. We know that $E$ is an extension of
all these stable factors. Therefore, it will be enough to show that whenever
there is an exact sequence $0 \to A \to B \to C \to 0$ with $\nu^{\min}_{\alpha,
\beta}(A), \nu^{\min}_{\alpha, \beta}(C) \geq a$, then
$\nu^{\min}_{\alpha, \beta}(B) \geq a$ for any $a \in \bbR$.

Assume there is a semistable quotient $B \onto D$ such that $\nu_{\alpha,
\beta}(D) < a$. Due to $\nu^{\min}_{\alpha, \beta}(A) > \nu_{\alpha, \beta}(D)$
there is no morphism from $A$ to $D$. Therefore, $B \onto D$ factors via a map
$C \to D$. But there is also no non trivial map from $C$ to $D$ because of
$\nu^{\min}_{\alpha, \beta}(C) > \nu_{\alpha, \beta}(D)$. But then $B \onto D$
is trival which is a contradiction.
\end{proof}
\end{lem}

This technical lemma allows to proceed with the proof of Theorem \ref{thm-main}.

\begin{lem}
The inclusion $\calC \subset \langle \calA_{\alpha,\beta}, \calA_{\alpha,\beta}[1]
\rangle$ holds.
\begin{proof}
If $L[i] \in \calB^{\alpha,\beta}$ holds for a line bundle $L$ and $i \in
\{0,1\}$, then $L[i]$ is tilt-stable (see Proposition 7.4.1 in \cite{BMT14}). By
Lemma \ref{computation} we get immediately $\calO(-1)[3], \calO [1],
\calO(1) \in \langle \calA_{\alpha,\beta}, \calA_{\alpha,\beta}[1] \rangle$.

By \cite{Ott88} the spinor bundle $S$ is $\mu$-stable. Since $\mu$-stability is
preserved by the tensor product (see \cite[Theorem 3.1.4]{HL10}) we obtain
$\mu$-stability of $S(-1)$. The inequality $\mu_{\alpha,\beta}(S(-1)) \leq 0$
leads to $S(-1)[1] \in \calB^{\alpha,\beta}$. In order to show $S(-1)[1] \in
\calA_{\alpha,\beta}$ we need to prove that any quotient $S(-1)[1] \onto G$ in
$\calB^{\alpha,\beta}$ satisfies $\nu_{\alpha,\beta}(G) > 0$. The proof proceeds
in three steps. At first we show $S(-1)[1]$ has no proper subobject of slope
$\infty$. Then we prove stability of $S(-1)[1]$ for $\beta=0$. Finally, we use
the previous lemma to reduce to this case.

Assume we have a proper subobject $A \into S(-1)[1]$ with $\nu_{\alpha,
\beta}(A) = \infty$. That means $\ch_1^{\beta}(A) = 0$ and moreover $\ch_1^{\beta}
(H^{-1}(A)) = 0$. Suppose we have $H^{-1}(A) \neq 0$. Then the injective
morphism $H^{-1}(A) \into S(-1)$ in $\Coh(Q)$ constitutes a contradiction to the
$\mu_{\alpha, \beta}$-stability of $S(-1)$ with the inequality $\mu_{\alpha,
\beta}(S(-1)) \leq 0$. Hence, $H^{-1}(A) = 0$ and since $\ch_1^{\beta}(A)=0$, it
follows that $A$ has rank $0$ and is supported in dimension less than or equal to one.
But in that case Serre duality implies $\Hom(A,S(-1)[1])=0$ which is a
contradiction to $A \to S(-1)[1]$ being a monomorphism.

Assume we have an exact sequence $0 \to A \to S(-1)[1] \to G \to 0$ in
$\calB^{\alpha,0}$ with $\nu_{\alpha,0}(G) \leq \nu_{\alpha,0}(A)$. The
long exact sequence in cohomology implies that $G\simeq N[1]$ for $N \in
\Coh{Q}$. Since $\ch_1(S(-1)[1]) = H$ (see Lemma \ref{computation}) and
$\nu_{\alpha,0}(G) \neq \infty$, we obtain $\ch_1(G) = H$ and $\ch_1(A)=0$. But
then $\nu_{\alpha,0}(A) = \infty$, a case that we had already ruled out.

Assume there is $\alpha_0 \in (0,1/3)$ and $\beta_0 \in [-1/2,0)$ such
that the inequality $\nu^{\min}_{\alpha_0, \beta_0}(S(-1)[1]) \leq 0$ holds.
Since stability is an open property by \cite[Corollary 3.3.3]{BMT14} and
$S(-1)[1]$ is $\nu_{\alpha_0, 0}$-stable, we get
$$\beta_1 := \sup \{\beta \leq 0 : \nu^{\min}_{\alpha_0, \beta}(S(-1)[1]) \leq 0
\} < 0.$$

Let $S(-1)[1] \onto N$ be a semistable  quotient in $\calB^{\alpha_0, \beta_1}$
as in Lemma \ref{lem:open}. Assume $\nu_{\alpha_0, \beta_1}(N) > 0$ and let
$N_1, \ldots, N_l$ be the stable quotients in the Jordan-H\"older filtration of
$N$. In a neighborhood around $(\alpha_0, \beta_1)$ we have the inequality
$\nu^{\min}_{\alpha, \beta}(S(-1)[1]) \geq \min \nu_{\alpha,\beta}(N_i) > 0$,
which is a contradiction to the choice of $\beta_1$. Therefore, we know
$\nu_{\alpha_0, \beta_1}(N) \leq 0$. We define the function
\begin{align*}
f(\beta) &= \frac{\alpha_0^2 H^2 \ch_1^{\beta_1}(N) \nu_{\alpha_0,
\beta}(N)}{\alpha_0}
\\
&= H \ch_2(N) - \beta H^2 \ch_1(N) + \frac{\beta^2 H^3}{2} \ch_0(N) -
\frac{\alpha_0^2 H^3}{2} \ch_0(N).
\end{align*}
We have the inequalities $f(\beta_1) \leq 0$ and $f'(\beta_1) = -H^2
\ch_1^{\beta_1}(N) < 0$. As $\nu^{\min}_{\alpha_0, \beta}(S(-1)[1]) \leq
\nu_{\alpha_0, \beta}(N)$ in a neighborhood of $(\alpha_0, \beta_1)$, the fact
that $f$ is decreasing at $\beta_1$ is a contradiction to the choice of
$\beta_1$.
\end{proof}
%\centerline{
%\begin{tikzpicture}
%  \draw plot[variable=\t,samples=1000, domain=0:72] ({sec(\t)},{tan(\t)})
%  node[above] {$C$};
%  \draw plot[variable=\t,samples=1000, domain=0:4] (\t,0);
%  \draw plot[variable=\t,samples=1000, domain=0:3] (2.5,\t);
%  \draw (1.63,0.3) node[anchor=south] {$\bullet^{(\alpha_0, \beta_1)}$};
%  \draw (2.919,0.3) node[anchor=south] {$\bullet^{(\alpha_0, 0)}$};
%\end{tikzpicture}
%}
%The set $C$ of all $(\alpha,\beta)$ with
%$\nu_{\alpha,\beta}(N[1]) = 0$ is a hyperbola given by the equation
%$$ f(\alpha,\beta) := \left( \beta -
%\frac{H^2\ch_1(N[1])}{H^3\ch_0(N[1])}\right)^2 - \alpha^2 -
%\frac{\overline{\Delta}_H(N[1])}{(H^3 \ch_0(N[1]))^2},$$ where
%$$\overline{\Delta}_H := (H^2 \ch_1)^2 - 2(H^3 \ch_0) (H\ch_2).$$
%By \cite[Corollary 7.3.2]{BMT14} we have
%$$\overline{\Delta}_H(N[1]) \geq 0.$$
%Therefore, $C$ has two branches left and right of the line $\beta
%= \frac{H^2\ch_1(N[1])}{H^3\ch_0(N[1])}$. But since $\ch_1^{\beta}(N[1]) > 0$
%holds, we only need to deal with the situation $\beta >
%\frac{H^2\ch_1(N[1])}{H^3\ch_0(N[1])}$. The inequality
% $\nu_{\alpha,\beta}(N[1]) \leq 0$ is equivalent to $f(\alpha,\beta) \geq 0$ as can be seen by a direct
%computation using $\ch_0(N[1]) \leq 0$. Therefore, $f(\alpha_0,\beta) \geq 0$
%for $\beta \geq \beta_1$ implies $\nu_{\alpha_0,\beta}(N[1]) < 0$. But
%that is a contradiction to the choice of $\beta_1$.
%\end{proof}
\end{lem}

The proof of Theorem \ref{thm-main} can be concluded by the next lemma.

\begin{lem}
For all $x \in X$, we have $k(x) \in \calC$ and for all proper subobjects $C
\into k(x)$ in $\calC$ the inequality $\Im Z_{\alpha,\beta,\frac{1}{6}}(C) > 0$
holds.
\begin{proof}
We have $k(x) \in \calC$ because of the resolution in \eqref{eq:resolution}
$$0 \to \calO(-1) \to S(-1)^{\oplus 2} \to \calO^{\oplus 4} \to \calO(1) \to
k(x) \to 0.$$
For the second assertion we need to figure out which are the subobjects of $k(x)
\in \calC$. Any object in $\calC$ is given by a complex $F$ of the form 
$$0 \to \calO(-1)^{\oplus a} \to S(-1)^{\oplus b} \to \calO^{\oplus c} \to
\calO(1)^{\oplus d} \to 0.$$
for $a,b,c,d \in \bbZ_{\geq 0}$. Since $\calC$ is the category of
representations of a quiver with relations with simple objects $\calO(-1)[3],
S(-1)[2], \calO[1], \calO(1)$, we can interpret $v(F)=(a,b,c,d)$ as the
dimension vector of that representation. Therefore, $F \into k(x) \onto G$
implies $a\leq 1$, $b\leq2$, $c\leq4$ and $d\leq 1$. If $F$ is non trivial, then
there is a simple object $T_1 \into F$. But the only simple object with non
trivial morphism into $k(x)$ is $\mathcal{O}(1)$. Therefore, the equality $d=1$
holds. If $G$ is non trivial, then there exists a simple quotient $k(x) \onto
T_2$. By Serre duality, the only simple quotient is $T_2 = \calO(-1)[3]$.
That implies $a=0$. Assume $b=2$, but $c<4$. Then we obtain $G = \calO(-1)[3]
\oplus \calO[1]^{\oplus 4-c} \onto \calO[1]$. A contradiction comes from
$\Hom(k(x), \calO[1]) = 0$. Therefore, $b=2$ implies $c=4$. The remaining cases
are $$v(F) \in \{(0,2,4,1)\} \cup \{(0,b,c,1) : b \in \{0,1\}, c \in
\{0,1,2,3,4\}\}.$$
Since $\Im Z_{\alpha,\beta,\frac{1}{6}}(S(-1))<0$, the case $b=0$ will follow
from $b=1$. With the same argument $v(F)=(0,1,4,1)$ will follow from
$v(F)=(0,2,4,1)$. Depending on the sign of $\Im
Z_{\alpha,\beta,\frac{1}{6}}(\calO)$, we can reduce the situation with $b=1$ to
either $c=0$ or $c=4$. Hence, we are left to check two cases. \\
\\
\centerline{
\begin{tabular}{ r | l }
  $v(F)$ & $\Im Z_{\alpha,\beta,\frac{1}{6}}$ \\
  \hline
  $(0,2,4,1)$ & $\alpha ((1 + \beta)^2 - \alpha^2)$ \\
  \hline
  $(0,1,0,1)$ & $\alpha (1 - 3(\beta^2 - \alpha^2))$
\end{tabular}} \\
\\
For all of them $\Im Z_{\alpha,\beta,\frac{1}{6}}$ is positive.
\end{proof}
\end{lem}

%The only remaining case is $\beta = -\alpha$. The issue in this case is that
%$Z_{\alpha,\beta,\frac{1}{6}}(\calO) = 0$. However, for all $0<\varepsilon
%\ll 1$ the previous proof still works for $s_0 = \frac{1}{6} + \varepsilon$.
%That is good enough.

\renewcommand{\refname}{References}

\end{document}